\documentclass[num-refs]{wiley-article}

\newcommand{\del}{\partial}
\def\MM#1{\boldsymbol{#1}}
\DeclareMathOperator{\diff}{d\!}

\usepackage{amscd}
\newcommand{\jump}[1]{\left[\!\left[ #1 \right]\!\right]}

\usepackage{color}
\usepackage{ulem}

\newcommand{\add}[1]{{#1}}

\usepackage{siunitx}

\papertype{Original Article}
\paperfield{Journal Section}

\title{Time parallel integration and phase averaging for the nonlinear
  shallow water equations on the sphere}

\author[1]{Hiroe Yamazaki}
\author[1]{Colin J. Cotter}
\author[2]{Beth A. Wingate}

\corraddress{Hiroe Yamazaki, Department of Mathematics, Imperial College London, London SW7 2AZ, UK}
\corremail{h.yamazaki@imperial.ac.uk}

\affil[1]{Department of Mathematics, Imperial College London, UK}
\affil[2]{Department of Mathematics, University of Exeter, UK}

\begin{document}
\maketitle

\begin{abstract}
We describe a proof-of-concept development and application of a phase
averaging technique to the nonlinear rotating shallow water equations
on the sphere, discretised using compatible finite element
methods. Phase averaging consists of averaging the nonlinearity over
phase shifts in the exponential of the linear wave operator. Phase
averaging aims to capture the slow dynamics in a solution that is
smoother in time (in transformed variables) so that larger timesteps
may be taken. We overcome the two key technical challenges that stand
in the way of studying the phase averaging and advancing its
implementation: 1) we have developed a stable matrix exponential
specific to finite elements and 2) we have developed a parallel finite
averaging proceedure. Following Peddle et al (2019), we consider
finite width phase averaging windows, since the equations have a
finite timescale separation. In our numerical implementation, the
averaging integral is replaced by a Riemann sum, where each term can
be evaluated in parallel. This creates an opportunity for parallelism
in the timestepping method, which we use here to compute our
solutions.  Here, we focus on the stability and accuracy of the
numerical solution. We confirm there is an optimal averaging window,
in agreement with theory. Critically, we observe that the combined
time discretisation and averaging error is much \add{smaller than} the time
discretisation error in a semi-implicit method applied to the same
spatial discretisation. An evaluation of the parallel aspects will
follow in later work.
\end{abstract}
\section{Introduction}

Phase averaging is a technique for approximating highly oscillatory
PDEs, such as the {{equations that govern the dynamics}} of models of
large scale geophysical fluid dynamics.  {{Examples of phase averaging
    include the solution of ODEs \citep{Sanders_etal_2007}, the
    analysis of fast singular limits \citep{majda1998averaging,
      KlMa1981} in simple geometries, and is an important concept used
    in understanding the influence of oscillations on mean flows,
    \add{with connections to Lagrangian averaging}. \add{In fact,
      theoretical and numerical work on phase averaging suggest that
      nonlinear phase averaging provides more accurate leading order
      dynamics,
      e.g. \citep{wagner_young_2015,kafiabad_vanneste_young_2021}. Therefore,
      the tools developed in this paper could contribute to further
      understanding this challenging topic in planetary fluid
      dynamics.} \add{In addition,} phase averaging has been proposed
    as a way to introduce more parallelization into numerical models
    \citep{haut2014asymptotic, Peddle2019PararealCF, Ariel_etal_16}.
    However, phase averaging in spherical geometries required for
    large scale geophysical fluid dynamics applications has two key
    technical barriers that stand in the way of both studying the
    phase average as a potentially interesting physical quantity as
    well as a potential basis for advancing parallel computing: 1) the
    development of a stable matrix exponential required for the
    mapping and 2) the development of a parallel finite interval
    averaging procedure. This paper overcomes these two challenges
    providing a proof-of-concept model of phase averaged time-stepping
    in spherical geometries, as well as a study of the method's
    accuracy.}}
    
We focus on the rotating shallow water equations in the low
Rossby number regime. The phase averaging technique averages the
nonlinearity over all phases of the fast waves to obtain an
approximation of the slow dynamics with no fast oscillations
present. Since the magnitude of time discretisation errors are
governed by the magnitude of time derivatives in the equation, this
means that these errors can be reduced for phase averaged models,
allowing larger timesteps can be taken in an accurate numerical
integration. However, there is an additional error introduced by the
phase averaging itself, which we shall investigate in this paper.

For the rotating shallow water equations, \citet{majda1998averaging}
showed that taking the low Rossby number limit in the phase averaged
equations leads to the quasigeostrophic equations.  However, the
quasigeostrophic equations are not uniformly valid, which is one
reason why their 3D counterpart is not used for operational weather
forecasting. \citet{haut2014asymptotic} proposed to use a phase
average of the nonlinearity over a finite width averaging window
$T$. For small $T$, the original equations are recovered and for large
$T$ the full phase averaging is recovered which filters all fast
dynamics. They proposed to perform the phase averaging numerically,
replacing the phase integral with a numerical quadrature rule. The
idea is to use parallel computation to implement the averaging: each
term in the quadrature rule can be evaluated independently and hence
in parallel. \citet{Peddle2019PararealCF} showed that, given a chosen
timestepping scheme and timestep size, there is an optimum averaging
window $T$ that minimises the total error (combining numerical
timestepping error and averaging error). Below the optimum, the
timestepping error dominates (and the timestepping scheme may become
unstable), and above the optimum, the averaging error dominates. For
larger timesteps, the optimal $T$ is larger, requiring more quadrature
points in the average and consequently requiring more computational
cores.

The finite window phase averaged model trades computational cores for
larger timesteps at the expense of accuracy (because of the additional
error from phase averaging). If the level of accuracy is insufficient,
a time-parallel predictor-corrector approach might be used to increase
accuracy. \citet{haut2014asymptotic} proposed to use the averaged
model as the coarse propagator in a highly efficient Parareal
iteration, demonstrating parallel speedups of a factor of 100 in a
rotating shallow water test case.  \citet{Peddle2019PararealCF} proved
convergence of the iterative Parareal procedure for highly oscillatory
PDEs with quadratic nonlinearity, making use of the optimal averaging
window $T$ at finite Rossby number. \citet{bauer2021higher} created a
hierarchy of higher order averaged models that increase accuracy
through increasing the number of variables; this type of hierarchy is
ideal for predictor-corrector approaches such as RIDC methods
\citep{ong2016algorithm} and PFASST \citep{minion2011hybrid} that
compute more accurate correction steps in parallel as new predictor
steps are being taken. These are motivations for the work in this
paper but here we focus on the impact of the averaging on the rotating
shallow water solution.

This paper addresses the challenge of producing a proof of concept of
the numerical averaging technique applied to the rotating shallow
water equations on the sphere. This required us to find a performant
way of implementing the necessary matrix exponentials, and to find a
discretisation approach that avoids instabilities in those
exponentials from spurious eigenvalues. It also required us to produce
a parallel implementation of the numerical averaging.  Our proof of
concept allows us to evaluate, for the first time, the
solution quality arising from the averaging technique against a test
case that is used in dynamical core development.

The rest of this paper is organised as follows.  In
Section \ref{sec:model}, we describe the phase averaging procedure,
and how it can be applied to the rotating shallow water equations. We
also describe our approach to timestepping these equations.  In
Section \ref{sec:results}, we present our numerical results,
examining the impact of averaging window $T$ and timestep $\Delta t$
on the errors associated with time integration. Finally,
in \ref{sec:conclusion} we provide a summary.

\section{Description of the method}
\label{sec:model}
\subsection{Shallow water equations}
In this section, we describe the shallow water equations and link them
to the general notation framework for phase averaging that will be
used in subsequent sections.

We begin with the nonlinear shallow water equations on a two
dimensional surface that is embedded in three dimensions,
\begin{eqnarray}
\MM{u}_t + f\MM{u}^{\perp} + (\MM{u}\cdot\nabla)\MM{u} + g\nabla\eta 
&=& 0, \label{ueq}\\
\eta_t + H\nabla\cdot\MM{u} + \nabla\add{\cdot}[\MM{u}(\eta - b)]
&=& 0, \label{etaeq}
\end{eqnarray}
where $\MM{u}$ is the horizontal velocity, $f$ is the Coriolis parameter, $\MM{u}^{\perp} = \MM{k} \times \MM{u}$ where $\MM{k}$ is the normal to the surface, and $g$ is the gravitational acceleration;
$\eta$ is the free surface elevation, $H$ is the mean layer thickness and $b$ is the height of the lower boundary, where the layer depth $h = H+\eta -b$; $\nabla$ and $\nabla\cdot$ are appropriate invariant gradient and divergence operators defined on the surface. Here we will concentrate on the case of the equations being solved on the surface of the sphere, so there are no boundary conditions to
consider.

Then we rewrite the equations as
\begin{eqnarray}
\MM{U}_t &=& \mathcal{L}\MM{U}+\mathcal{N}(\MM{U}), \label{split_eqs}
\end{eqnarray}
where the vector of unknowns $\MM{U}(t) = (\MM{u}, \eta)$. The matrix $\mathcal{L}$ represents a linear operator and $\mathcal{N}(\cdot)$ is a nonlinear operator which satisfy
\begin{eqnarray}
\mathcal{L}\MM{U} 
&=& \left(\begin{array}{cc} -f(\cdot)^{\perp} & -g\nabla\\ -H\nabla\cdot & 0\\ \end{array}\right) 
\left(\begin{array}{c} \MM{u}\\ \eta\\ \end{array}\right)
= \left(\begin{array}{cc} -f(\MM{u})^{\perp} & -g\nabla\eta\\ -H\nabla\cdot\MM{u} & 0\\ \end{array}\right),
\label{linear_operator}\\
\mathcal{N}(\MM{U})
&=& \left(\begin{array}{c} -(\MM{u}\cdot\nabla)\MM{u}\\ -\nabla\add{\cdot}\left[\MM{u}(\eta - b)\right]\\ \end{array}\right).
\label{nonlinear_operator}
\end{eqnarray}

\subsection{Phase averaging}
Now we consider an approximation to the equation \eqref{split_eqs} by averaging the nonlinearity over the fast oscillations.
First we introduce a coordinate transformation,
\begin{eqnarray}
\MM{V}(t)&=&e^{-\mathcal{L}t}\MM{U}(t), \label{v_eq}\\
\frac{\del\MM{V}}{\del t}(t) &=& e^{-\mathcal{L}t}\mathcal{N}
\left(e^{\mathcal{L}t}\MM{V}(t)\right), \label{vt_eq}
\end{eqnarray}
where {$e^{\mathcal{L}t}\MM{V}$ is solution at time $t$ to the linear part of the equation $\MM{U}_t = \mathcal{L}\MM{U}$ with an initial condition of $\MM{U}(0) = \MM{V}$.}

To allow averaging the model over a finite time interval $s$ around
time $t$, we follow \citet{bauer2021higher} and extend the equations
\eqref{v_eq} and \eqref{vt_eq} by introducing a phase variable $s$,
\begin{eqnarray}
\MM{V}(t,s)&=&e^{-\mathcal{L}(t+s)}\MM{U}(t, s), \label{v_eq_s}\\
\frac{\del\MM{V}}{\del t}(t, s) &=& e^{-\mathcal{L}(t+s)}\mathcal{N}
\left(e^{\mathcal{L}(t+s)}\MM{V}(t, s)\right). \label{vt_eq_s}
\end{eqnarray}
\add{Of course, we are only interested in the solution at $s=0$, but
  the averaging approximation below corresponds to the assumption that the
  solution $V(t,s)$ is insensitive to changes in $s$. \citet{bauer2021higher}
  introduced higher order approximations that parameterise the sensitivity
  to $s$ near $s=0$.}
  
An averaging approximation to \eqref{vt_eq_s} over the averaging window $T$ with respect to a weight function $\rho$ can be written as
\begin{eqnarray}
\frac{\del\bar{\MM{V}}}{\del t}(t) &=& \add{\frac{1}{2T}}\int^{T}_{-T}\rho\left(\frac{s}{T}\right)e^{-\mathcal{L}(t+s)}\mathcal{N}\left(e^{\mathcal{L}(t+s)}\bar{\MM{V}}(t)\right) \diff s, \label{vt_averaged}
\end{eqnarray}
where $\bar{\MM{V}}$ denotes the averaged $\MM{V}$. In our computations, $\rho$ is given by
\add{
  \begin{equation}
    \rho(s) = \left\{
    \begin{array}{rl}
\rho_0\exp(-1/((s-0.5)(s+0.5))), & \quad -0.5 < s < 0.5, \\
0, & \quad |s| \geq 0.5, \\
\end{array}\right.
\end{equation}
with an appropriate weighting factor,} but other weight functions may
be used provided that they integrate to 1.

\add{It is very important to note that this averaging integral is \emph{not}
  integrating along the history of $\bar{V}$ (which would be an integral in
  the time variable $t$), but is instead integrating over the phase shift
  variable $s$ in the exponential operators. Equation \eqref{vt_averaged}}
  is similar to the phase averaging in \cite{Peddle2019PararealCF}. The main difference is that our phase shift is defined in the mapping \eqref{v_eq_s}, whereas \cite{Peddle2019PararealCF} introduce the phase shift parameter, $s$, in the nonlinear term.

\add{The unaveraged equation \eqref{vt_eq_s} has oscillatory
  nonautonomous contributions coming from the exponentials
  (corresponding to fast waves in the untransformed equations
  \eqref{split_eqs}), and the averaging approximation filters out the
  components of these contributions with time period below $T$.  This
  approximation is best understood in various limits. In the limit as
  $T\to 0$, $\rho(s/T)/(2T) \to \delta(s)$, so we recover the original
  unaveraged equation \eqref{vt_eq}. In the limit $T\to \infty$, the
  asymptotic approximation $U(t) =
  \exp(-t/\epsilon\mathcal{L}_0)\bar{V}(t) + \mathcal{O}(\epsilon)$ is
  obtained, where $\mathcal{L}=\epsilon \mathcal{L}_0$
  \citep{majda1998averaging}. This asymptotic approximation describes
  a slowly evolving solution $\bar{V}(t)$ with superposed fast linear
  waves: the classic ``slow manifold'' picture. For intermediate
  values of $T$, we can select which nonautonmous contributions we
  want to remove and which contributions we want to retain. This can
  be important if the interaction between the two exponentials inside
  and outside $\mathcal{N}$ in \eqref{vt_eq_s} (\emph{e.g.} through
  triad interactions in the case of quadratic nonlinearity) leads
  to oscillations with frequencies approaching those of the dynamics
  of $\bar{V}$; in that case we would lose long time accuracy if we increased
  $T$ to remove those oscillations.}

For an implementable method, we replace the integral with a
Riemann sum to obtain
\begin{eqnarray}
\frac{\del\bar{\MM{V}}}{\del t} \!\!\!\!
&\,& \,\simeq \quad {\sum_{k=-N}^{N}} w_{k}e^{-\mathcal{L}(t+s_{k})}\mathcal{N}
\left(e^{\mathcal{L}(t+s_{k})}\bar{\MM{V}}(t)\right) \\
&\,& := \quad \left\langle e^{-\mathcal{L}(t+s)}\mathcal{N}\left(e^{\mathcal{L}(t+s)}\bar{\MM{V}}(t)\right)\right\rangle_s,  \label{vt_discretised}
\end{eqnarray}
where $w_k$ are appropriate weight coefficients (obtained from the
product quadrature rule weights and the value of the weight function
$\rho$ evaluated at the quadrature points) and $s_k = kT/N$. This
defines the angle bracket notation \eqref{vt_discretised}, used in the
following section.  \add{Since this is an oscillatory integral, there
  is little reduction in the quadrature error until the oscillations
  are resolved on the quadrature points, after which the error
  collapses quickly. In this work we used equispaced quadrature
  points, with four points per time period of the fastest frequency of
  $\mathcal{L}$; the resulting dynamics was very insensitive to
  increasing this number beyond four.}

\add{Whilst \eqref{vt_discretised} looks complicated, each of the terms
  in the sum is independent and so they can be evaluated in parallel
  if computational resources are available.
  }

\subsection{Time discretisation}
In this section we describe our time integration approach. The general
summary is that we use an averaged version of a Lawson exponential
integrator (see \citet{hochbruck2010exponential} for a review). This
means that we apply a standard time integration method to
\eqref{vt_discretised}, and then transform from $\MM{V}$ back to
$\MM{U}$ to restrict exponentiation to time intervals of
$\mathcal{O}(\Delta t)$. These shorter exponentiations are less
expensive to compute numerically.

For the classical 4th order Runge-Kutta
scheme, we obtain
\begin{align}
\bar{\MM{U}}_{1} &=& \bar{\MM{U}}^{n} + \frac{\Delta t}{2} \, \left\langle e^{-\mathcal{L}s}\mathcal{N}(e^{\mathcal{L}s}\bar{\MM{U}}^{n})\right\rangle_s, 
\label{stage1_c}\\
\bar{\MM{U}}_{2}
&=& e^{\mathcal{L}\frac{\Delta t}{2}}\bar{\MM{U}}^{n} + \frac{\Delta t}{2} \left\langle e^{-\mathcal{L}s}\mathcal{N}(e^{\mathcal{L}\frac{\Delta t}{2}}e^{\mathcal{L}s}\bar{\MM{U}}_1)\right\rangle_s, 
\label{stage2s}\\
\bar{\MM{U}}_{3}
&=& e^{\mathcal{L}\frac{\Delta t}{2}}\bar{\MM{U}}^{n} + \Delta t \left\langle e^{-\mathcal{L}s}\mathcal{N}(e^{\mathcal{L}s}\bar{\MM{U}}_2)\right\rangle_s, 
\label{stage3s}\\
\bar{\MM{U}}^{n+1}
&=& e^{\mathcal{L}\Delta t}\bar{\MM{U}}^{n} + \frac{\Delta t}{6}\Bigg[e^{\mathcal{L}\Delta t}\left\langle e^{-\mathcal{L}s}\mathcal{N}(e^{\mathcal{L}s}\bar{\MM{U}}^n)\right\rangle_s
+ 2 e^{\mathcal{L}\frac{\Delta t}{2}}\left\langle e^{-\mathcal{L}s}\mathcal{N}(e^{\mathcal{L}\frac{\Delta t}{2}}e^{\mathcal{L}s}\bar{\MM{U}}_1)\right\rangle_s \nonumber \\
& &+ 2 e^{\mathcal{L}\frac{\Delta t}{2}}\left\langle e^{-\mathcal{L}s}\mathcal{N}(e^{\mathcal{L}s}\bar{\MM{U}}_2)\right\rangle_s
+ \left\langle e^{-\mathcal{L}s}\mathcal{N}(e^{\mathcal{L}\frac{\Delta t}{2}}e^{\mathcal{L}s}\bar{\MM{U}}_3)\right\rangle_s\Bigg].
\label{stage4}
\end{align}

\subsection{Chebyshev exponentiation}
Implementing exponentials $e^{\mathcal{L}t}$ of grid based
discretisations on the sphere is challenging, because the efficiency
and parallel scalability of these discretisations relies upon matrix
sparsity, and $e^{\mathcal{L}t}$ is not sparse. Instead we need to
consider scalable algorithms that construct $e^{\mathcal{L}t}$ using
only sparse matrix applications and local operations. This is done for
the first time in this paper in the context of numerical averaging
techniques applied to the rotating shallow water equations on the
sphere.

To implement the exponential operator $e^{\mathcal{L}t}$, we use a
Chebyshev approximation,
\begin{eqnarray}
e^{\mathcal{L}t}\bar{\MM{U}} \approx \sum^{N}_{k=0}a_{k}P_{k}(\mathcal{L}t)\bar{\MM{U}},
\label{cheby}
\end{eqnarray}
where $N$ is the number of polynomials, $a_k$ are polynomial
coefficients, and $P_k$ are modified Chebyshev polynomials. The modification
is a change of coordinates transforming the imaginary interval
\begin{equation}
  \left\{ z = iy: y \in [-L, L]\right\},
\end{equation}
to the unit interval $[-1,1]$, where $L=|\lambda_{\max}|t_{\max}$,
$\lambda_{\max}$ is the eigenvalue of $\mathcal{L}$ with maximum magnitude, and
the approximation is valid for times $|t|<t_{\max}$.

\add{In general, Krylov subspace methods (of which the Chebyshev
  approach is one) for oscillatory problems require a number of
  iterations proportional to the Courant number
  \citet{hochbruck2010exponential,pieper2019exponential}.
The advantage with Chebyshev polynomials is that there is a three
term recurrence, so we do not need to store and compute with the
entire Krylov basis.}
\add{This method starts with the three term recurrence for
Chebyshev polynomials,} modified with the above coordinate
transformation to get
\begin{eqnarray}
P_0(l) &=& 1, \\
P_1(l) &=& \frac{-il}{L}, \\
P_{k+1}(l) &=& \frac{2l P_k(l)}{iL} - P_{k-1}(l).
\end{eqnarray}
Therefore $P_k(\mathcal{L}t)\bar{\MM{U}}$ in the equation
\eqref{cheby} are obtained recursively as
\begin{eqnarray}
P_0(\mathcal{L}t)\bar{\MM{U}} &=& \bar{\MM{U}}, \\
P_1(\mathcal{L}t)\bar{\MM{U}} &=& \frac{-it\mathcal{L}\bar{\MM{U}}}{L}, \\
P_{k+1}(\mathcal{L}t)\bar{\MM{U}} &=& \frac{2t P_k(\mathcal{L}t)\mathcal{L}\bar{\MM{U}}}{iL} - P_{k-1}(\mathcal{L}t)\bar{\MM{U}}.
\end{eqnarray}
This avoids explicitly forming polynomials of matrices, by instead
just recursively forming the action of polynomials of matrices on
vectors\add{, by repeated application of $\mathcal{L}$. Further,
  this application is performed matrix free by writing the action
  of the matrix on a vector equivalently as a set of integrals
  (in the usual finite element manner), although a mass matrix
  solve is required in this formulation.}
  
As usual for Chebyshev polynomials, the coefficients $a_k$
are computed by using a fast Fourier transform
\citep{trefethen2019approximation}.  Following
\citet{gander2013paraexp}, we then discard coefficients starting at
the highest degree and going downwards until the total magnitude of
discarded coefficients exceeds some threshold ($1.0e-6$, for our
results).

We note that for larger averaging windows, higher degree Chebyshev
approximations are needed, meaning that the application of the
approximated exponential operator takes longer. \add{This is illustrated
  in Table {\ref{tab:cheb}}.}

\begin{table}
  \add{
    \begin{center}
  \begin{tabular}{|c|c|c|c|c|}
    \hline
    & Mesh refinement level & $\Delta s=15$ minutes & $\Delta s=30$ minutes & $\Delta s=1$ hour \\
    \hline
    Number of iterations & 3 & 10 & 15 & 24 \\
    Number of iterations & 4 & 16 & 25 & 41 \\
    Number of iterations & 5 & 26 & 42 & 73 \\
    \hline
  \end{tabular}
  \end{center}}
  \caption{\label{tab:cheb}\add{Table showing the number of Chebyshev
      iterations required to provide a Chebyshev approximation of the
      exponential that contains all but a $1.0e-6$th of the energy
      of the full exponential, for different mesh refinements and
      times. We see that the growth is approximately linear
      in both mesh resolution and time.}}
\end{table}

\add{For a fully performant method, in the future we will incorporate}
the REXI technique of \citet{haut2016high}, which approximates the
exponential along the imaginary axis by a sum of rational
polynomials. This method can be parallelised over the sum; each
parallel term requires the solution of a complex valued elliptic
problem of the form $aI + b\mathcal{L}$, where $a$ and $b$ are
rational coefficients.  Some more details of the implementation and
examination of parallel performance are provided in
\citet{schreiber2018beyond}. When extending to three dimensions, it
may also be useful to exploit the vertical horizontal tensor product
structure in the exponential, as discussed in
\citet{croci2022exploiting}.

However, here we are focussed on the error behaviour of the averaging
technique, so the Chebyshev approximation suffices for this purpose.
\add{The goal here is to use parallelism as a way of understanding the
  impact of the averaging technique on the solution of the PDE.}

\subsection{Spatial discretisation}

In this study we used the compatible finite element discretisation for
the nonlinear rotating shallow water equations on the sphere given in
\citet{gibson2019compatible}. This was chosen because it leads to a
discretised $\mathcal{L}$ that still has purely imaginary eigenvalues,
and has a discrete Helmholtz decomposition that correctly separates
the fast inertia gravity waves and the slow balanced motion. This is
critical to addressing the challenge of applying numerical averaging
to the rotating shallow water equations on the sphere. Any similar
approach to the discretisation with these properties (e.g. a spectral
discretisation, or C-grid finite difference method) is expected to
produce similar results.

The compatible finite element discretisation is built around a pair of
spaces $\mathbb{V}_1\subset H(\mathrm{div})$ and $\mathbb{V}_2\subset
L_2$, selecting $\MM{u}\in \mathbb{V}_1$ and $\eta \in \mathbb{V}_2$.
In these examples we chose BDM2 for $\mathbb{V}_1$ and P1$_{\mathrm{DG}}$
for $\mathbb{V}_2$, producing a second-order scheme in space.

The discrete linear operator
$\mathcal{L}:\mathbb{V}_1\times\mathbb{V}_2 \to
\mathbb{V}_1\times\mathbb{V}_2$ is then defined by
$\mathcal{L}(\MM{u},\eta)=(\MM{u}_1,\eta_1)$, where
\begin{eqnarray}
\int_{\Omega} \MM{w} \cdot \MM{u}_1 \diff x &=& - \int_{\Omega} \MM{w} \cdot (f\MM{u}^{\perp}) \diff x + \int_{\Omega} (\nabla \cdot \MM{w})g \eta \diff x, \quad \forall\MM{w} \in \mathbb{V}_1, \label{ueq_int}\\
\int_{\Omega} \phi \, \eta_1 \diff x &=& - H \int_{\Omega} \phi \nabla \cdot \MM{u} \diff x, \quad \forall\phi \in \mathbb{V}_2.
\end{eqnarray}
Implementing this requires the solution of a block diagonal system for
the basis coefficients of $\eta_1$ and a sparse (but globally coupled)
system for the basis coefficients of $\MM{u}_1$. For the latter we
observe a mesh-independent number of iterations when solving using
a scalable iterative method (described below).

The discrete nonlinear operator 
$\mathcal{N}:\mathbb{V}_1\times\mathbb{V}_2 \to
\mathbb{V}_1\times\mathbb{V}_2$ is then similarly defined by
$\mathcal{N}(\MM{u},\eta)=(\MM{u}_2,\eta_2)$, where
\begin{eqnarray}
\int_{\Omega} \MM{w} \cdot \MM{u}_2 \diff x &=& \int_{\Omega} \MM{u} \cdot \nabla^{\perp} (\MM{u}^{\perp} \cdot \MM{w}) \diff x - \int_\Gamma \jump{\MM{n}^{\perp}(\MM{u}^{\perp}\cdot\MM{w})} \cdot \tilde{\MM{u}} \diff S + \int_{\Omega} \nabla \cdot \MM{w} \left(\frac{1}{2}|\MM{u}|^{2}\right) \diff x, \quad \forall\MM{w} \in \mathbb{V}_1,\\
\int_{\Omega} \phi \, \eta_2 \diff x &=& \int_{\Omega} \nabla \phi \, \MM{u} (\eta - b) \diff x + \int_{\Gamma} \jump{\phi\MM{u}}(\eta - b) \diff S, \quad \forall\phi \in \mathbb{V}_2,
\end{eqnarray}
where $\MM{n}$ is the outward pointing unit normal vector to the boundary $\del \Omega$ of $\Omega$, $\Gamma$ denotes the set of interior facets in the mesh with the two sides of each facet arbitrarily labeled by $+$ and $-$, the jump operator is defined by
\begin{eqnarray}
\jump{q} &=& q^{+}\MM{n}^{+} + q^{-}\MM{n}^{-}, \\
\jump{\MM{v}} &=& \MM{v}^{+}\cdot\MM{n}^{+} + \MM{v}^{-}\cdot\MM{n}^{-},
\end{eqnarray}
for any scalar $q$ and vector $\MM{v}$, and $\tilde{\MM{u}}$ is evaluated on the
upwind side as
\begin{eqnarray}
\tilde{\MM{u}} = \left\{
\begin{array}{l}
\MM{u}^{+}  \ \ \ \ \mbox{if} \ \MM{u}\cdot\MM{n}^{+} < 0, \\
\MM{u}^{-} \ \ \ \ \mbox{otherwise.}
\end{array}
\right.
\end{eqnarray}
Implementing this requires the solution of the same systems for
$\MM{u}_2$ and $\eta_2$ as $\MM{u}_1$ and $\eta_1$, respectively.
\add{This upwind stabilisation of the advection terms is the
  only dissipative term and there are no explicit dissipation
  terms in the model.}

Our code implementation was written using Firedrake
\citep{rathgeber2016firedrake}, an automated system for the solution
of partial differential equations using the finite element method,
with the resulting matrix systems being solved using PETSc
\citep{balay2021petsc}. A direct solver was used for the block
diagonal systems for $\eta$ and the conjugate gradient method
preconditioned by incomplete Cholesky factorisation was used to solve
the sparse systems for $\MM{u}$. The terms of the average are computed
in parallel using the ``ensemble parallelism'' capability of
Firedrake, which was implemented for this project. This implementation
provides MPI subcommunicators for the distribution of the terms of the
averaging sum, with the sum being formed by reduction over the
subcommunicators.

\section{Numerical experiments}\label{sec:results}

In this section we show numerical results from a standard test case on the sphere described by \citet{williamson1992standard}.
Here we use their test case number 5 (flow over a mountain), where the model is initialised with the layer depth and velocity fields that are in geostrophic balance: 
\begin{eqnarray}
h &=& H - \left(R\Omega u_0 + \frac{u_0^{2}}{2}\right)\frac{z^2}{gR^2},\\
\MM{u} &=& \frac{u_0}{R}(-y, x, 0),
\end{eqnarray}
where $R = 6.37122 \times 10^6$ m is the radius of the Earth, $\Omega = 7.292 \times 10^5$ s$^{-1}$ is the rotation rate of the Earth, $(x, y ,z)$ are the 3D Cartesian coordinates, the maximum zonal wind speed $u_0$ = 20 m, $g$ = 9.8 m s$^{-2}$ and $H$ = 5960 m. 
An isolated mountain is placed with its centre at latitude $\phi = \pi/6$ and longitude $\lambda = -\pi/2$. The height of the mountain is described as 
\begin{eqnarray}
b &=& b_0\left(1-\frac{(\mathrm{min}[R_0^2, (\phi-\phi_c)^2+(\lambda-\lambda_c)^2])^{1/2}}{R_0}\right)
\end{eqnarray}
where $b_0$ = 2000 m and $R_0 = \pi / 9$. The sudden appearance of
this mountain disturbing the balanced flow creates significant fast
unbalanced inertia-gravity waves as well as triggering slow balanced
vortex motion.

Icosahedral grids with a piecewise cubic approximation to the sphere
are used in the model.  The number of cells is \add{$N_C=20480$, and} the
maximum cell centre to cell centre distance is 263 km and the minimum
distance is 171 km.  A timestep of 900 s is used in the averaged model
unless stated otherwise.  As there is no analytical solution for this
problem, the model output is compared to a reference solution
generated from a semi-implicit nonlinear shallow water code provided
by \citet[][]{gibson2019compatible} (which we refer to as the
``standard model''), using the same spatial resolution of $N_C = 20480$,
\add{and the same spatial discretisation}.  A much smaller timestep of 22.5
s is used to generate the reference solution, to get as close as
possible to the exact solution to the time continuous space discrete
system we are trying to approximate with the averaging technique.

\begin{figure}
  \centering
  \includegraphics[height=40mm]{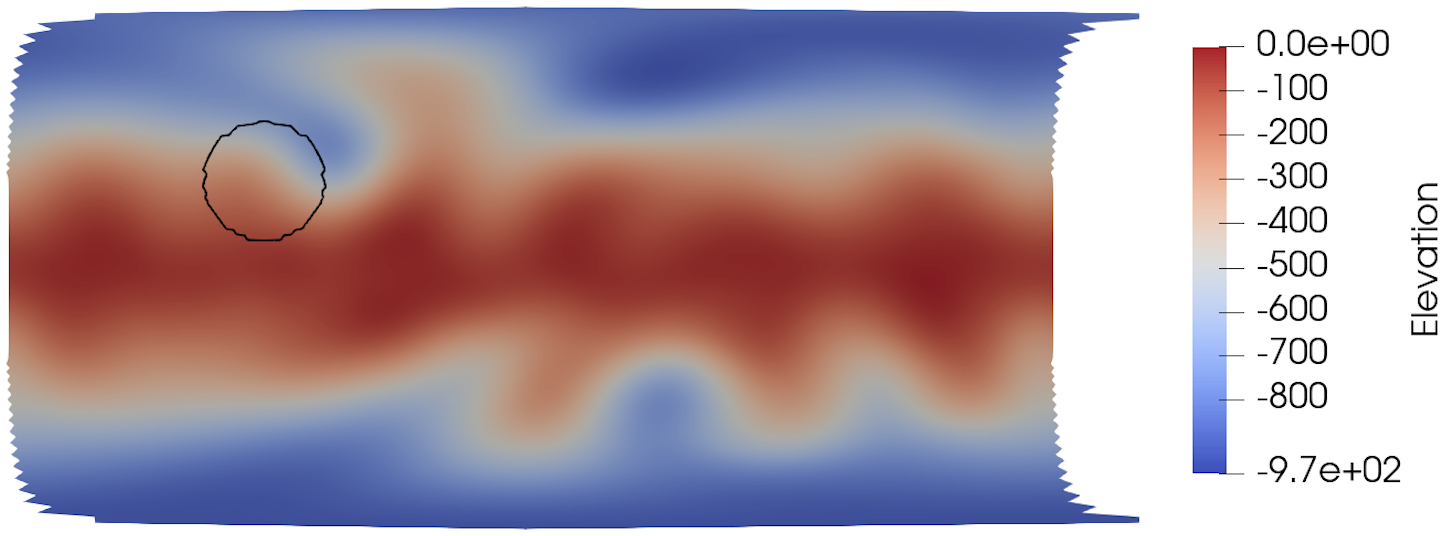}
  \caption{Free surface elevation $\eta$ at day 15 from the averaged model. The timestep is $\Delta t = 900$ s and the averaging window is $T = 1$ \add{hour}. The solid line indicates the position of the mountain.}
  \label{fig:eta}
\end{figure}

\begin{figure}
  \centering
  \includegraphics[height=40mm]{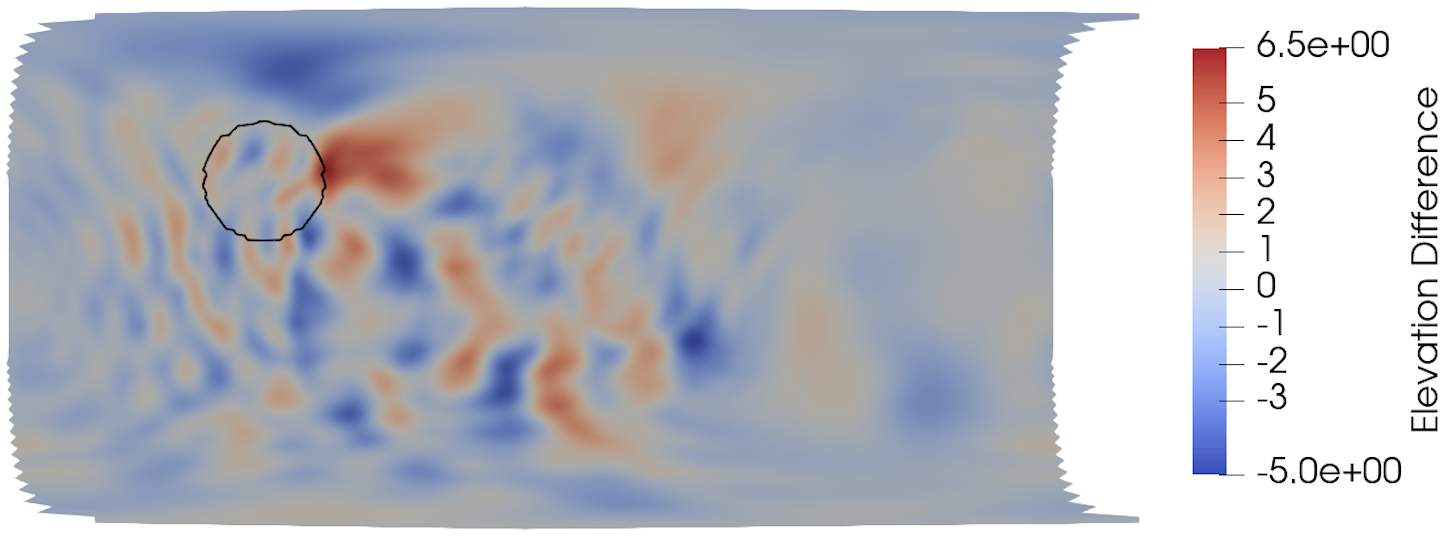}
  \caption{Errors in the elevation $\eta$ at day 15 compared to the reference solution. The solid line indicates the position of the mountain.}
  \label{fig:etadiff}
\end{figure}

Figure \ref{fig:eta} shows the field of the free surface elevation $\eta$ at day 15 from the averaged model plotted in a latitude longitude space.
The averaging window is $T$ = 1 hour in this plot.
The model successfully reproduces waves that travel around the globe as a result of the zonal flow interacting with the mountain.
Figure \ref{fig:etadiff} shows the errors in $\eta$ at day 15 compared to the reference solution.
We can see that errors are sufficiently small and not dominated by errors due to grid imprinting.

\begin{figure}
  \centering
  \includegraphics[height=40mm]{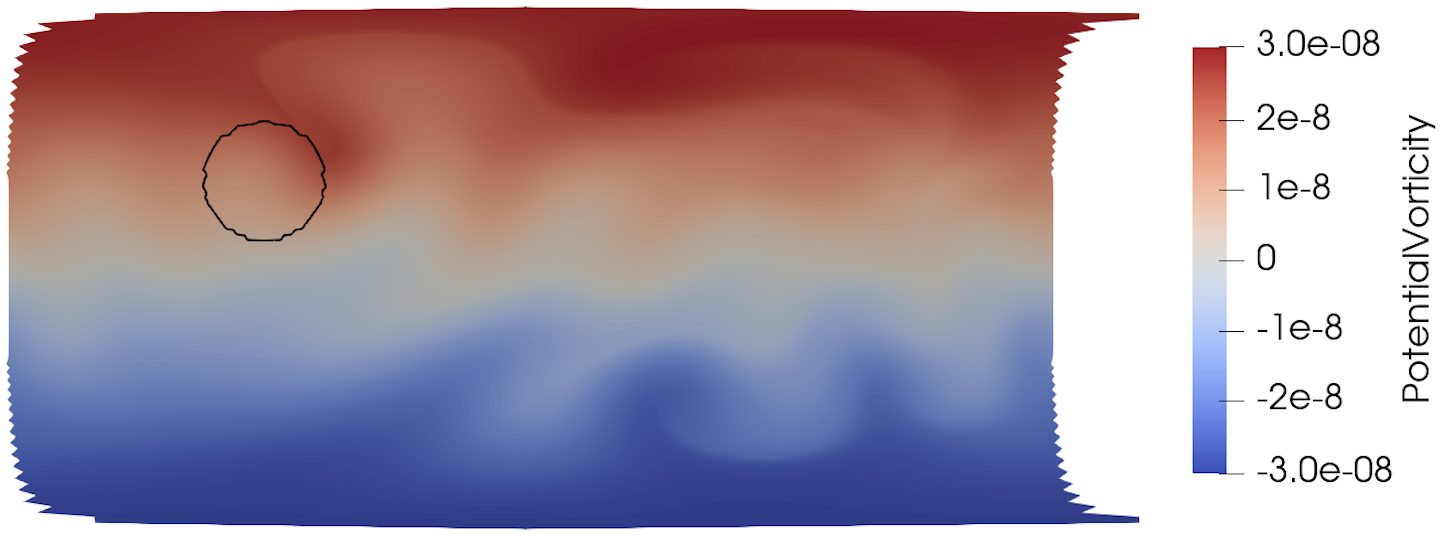} \\
  \vspace{0.5em}
  \includegraphics[height=40mm]{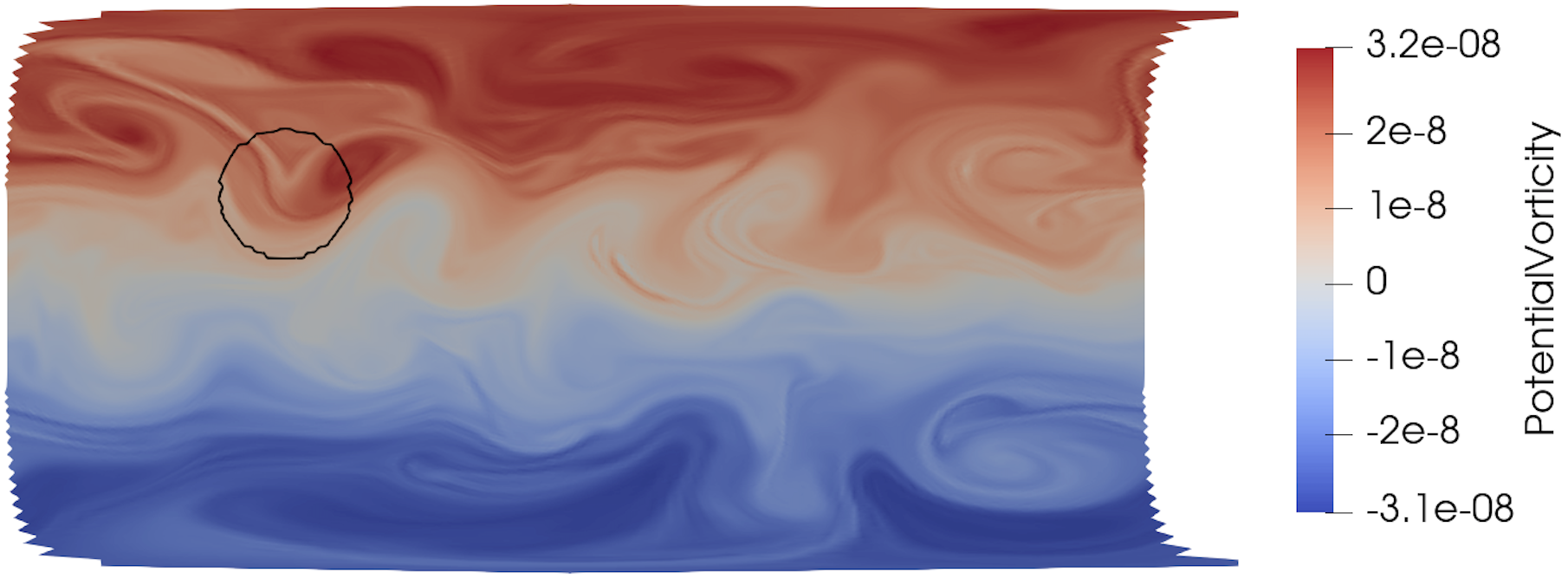}
  \caption{Potential vorticity at (top) day 15 and (bottom) day 50 from the averaged model. The timestep is $\Delta t = 900$ s and the averaging window is $T = 1$ hour. The solid lines indicate the position of the mountain.}
  \label{fig:pv}
\end{figure}

Figure \ref{fig:pv} shows the fields of the potential vorticity at day 15 and 50 from the averaged model. The flow is only weakly nonlinear at day 15, and fine scale structure has been generated at day 50 as the flow becomes more nonlinear. 
These results are consistent with the numerical results by \citet[][]{thuburn2014mimetic} and \citet[][]{shipton2018higher}.  

\begin{figure}
  \centering
    \includegraphics[height=70mm]{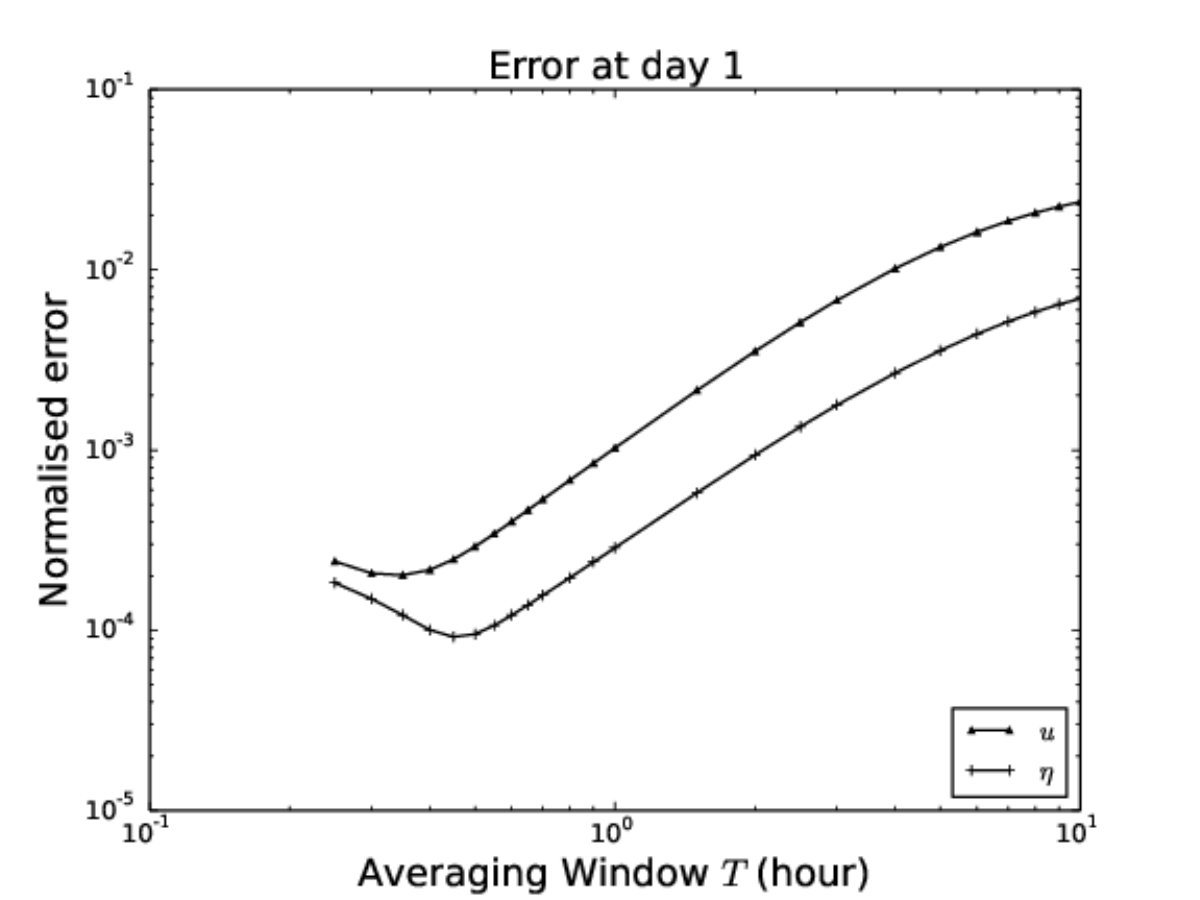}
  \caption{Normalised errors in the averaged model at day 1 versus the averaging window $T$. The solid {curve} with triangles shows the H(div) norm of $\MM{u}$, and the solid {curve} with crosses shows the $L_2$ norm of $\eta$, normalised by the norms of the reference solution. The time step $\Delta t = 900$ s is fixed in all the simulations. Note the clear existence of optimal averaging windows for each variable.
  \label{fig:peddle_plots}}
\end{figure}
\begin{figure}
    \includegraphics[height=70mm]{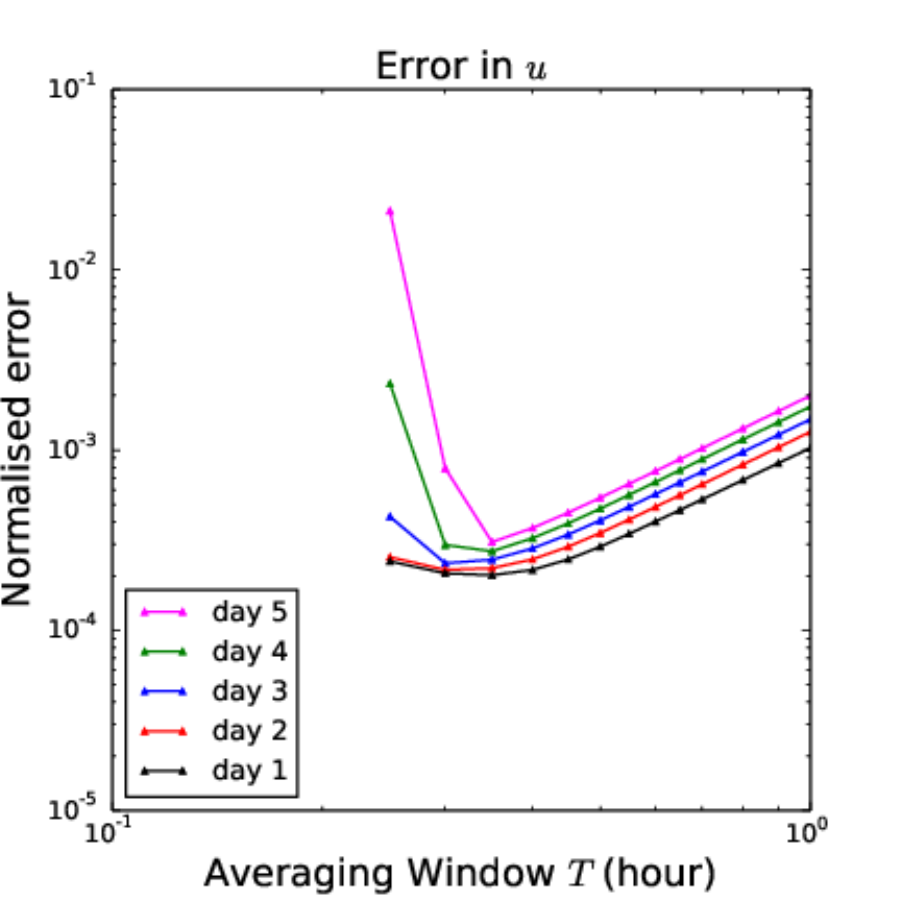} 
    \includegraphics[height=70mm]{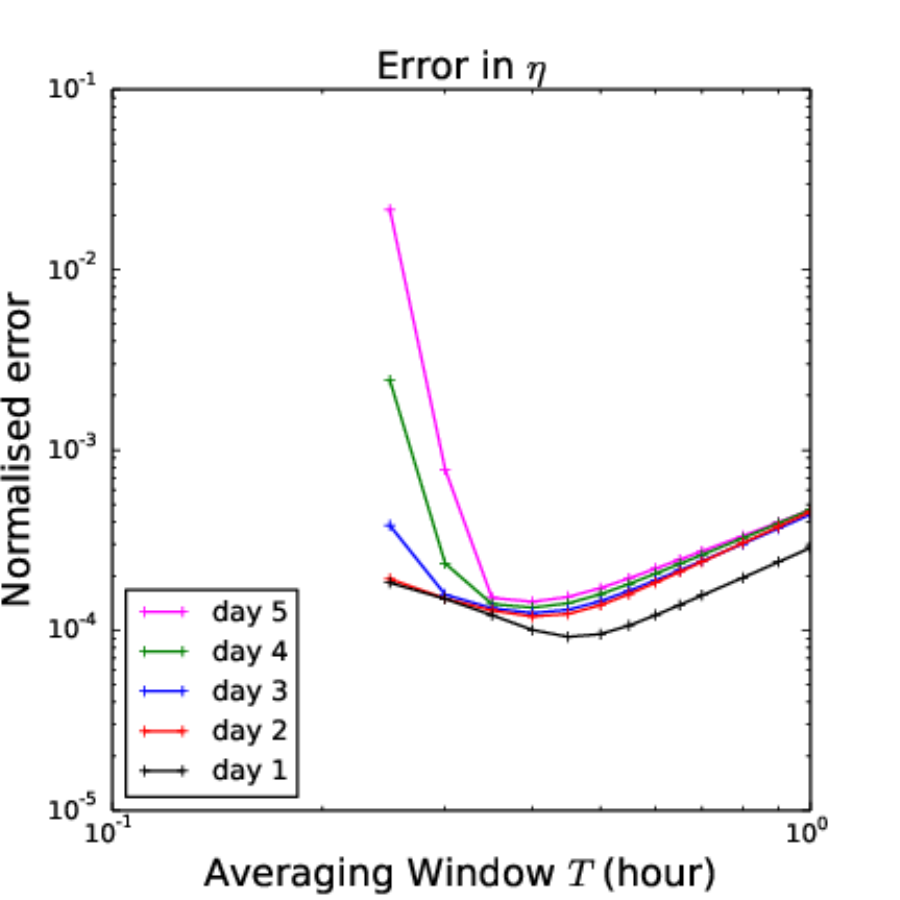} 
\caption{ Time evolution of the errors up to day 5: (left) the H(div)
  norm of $\MM{u}$ and (right) the $L_2$ norm of $\eta$, normalised by
  the norms of the reference solution. {The five curves} in each panel
  show the corresponding errors at day 1, 2, 3, 4 and 5 from the
  bottom to the top, respectively. The time step $\Delta t = 900$ s is
  fixed in all the simulations. 
  \label{fig:peddle_plots_5days}}
\end{figure}

Now, we examine the impact of the averaging window $T$ on the accuracy of the averaged model. Figure \ref{fig:peddle_plots} shows the H(div) norm of $\MM{u}$ and the $L_2$ norm of $\eta$ at day 1, normalised by the norms of the reference solution, plotted over various averaging windows.
The spatial resolution of $N = 20480$ and the time step of $\Delta t = 900$ s were kept the same as in Figure \ref{fig:etadiff}, whereas a range of values
between 0.25 hours and 10 hours were used for the the averaging window $T$.
The result reveals the clear existence of optimal averaging windows at around $T = 0.35$ hours for $\MM{u}$ and $T = 0.45$ hours for $\eta$, respectively; this will vary in the choice of norm.
This result demonstrates that the behavior in the averaged model is consistent with the error bounds shown in \citet[][]{Peddle2019PararealCF}.
Figure \ref{fig:peddle_plots_5days} shows the time evolutions of the errors in $\MM{u}$ and $\eta$ up to day 5.
For both variables, the minima in the error curves move to the left at day 2. From day 3, the errors at averaging windows $T$ $\leq$ 0.3 grow rapidly, showing that the model is slowly blowing up at those small averaging windows. This is because the fixed timestep is not resolving the unfiltered fast oscillations
at these small averaging window widths.
\begin{figure}
  \includegraphics[height=65mm]{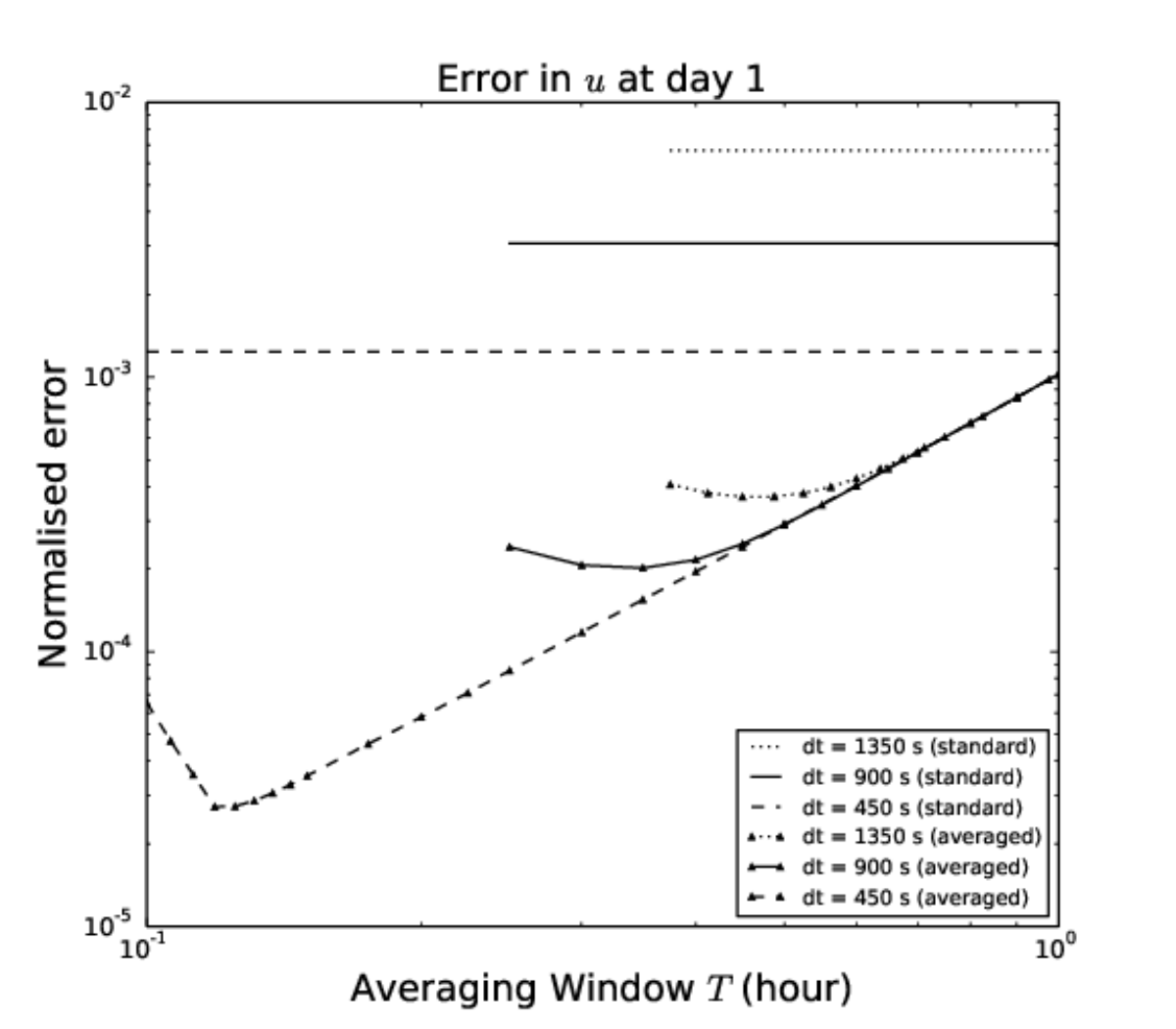} 
  \includegraphics[height=65mm]{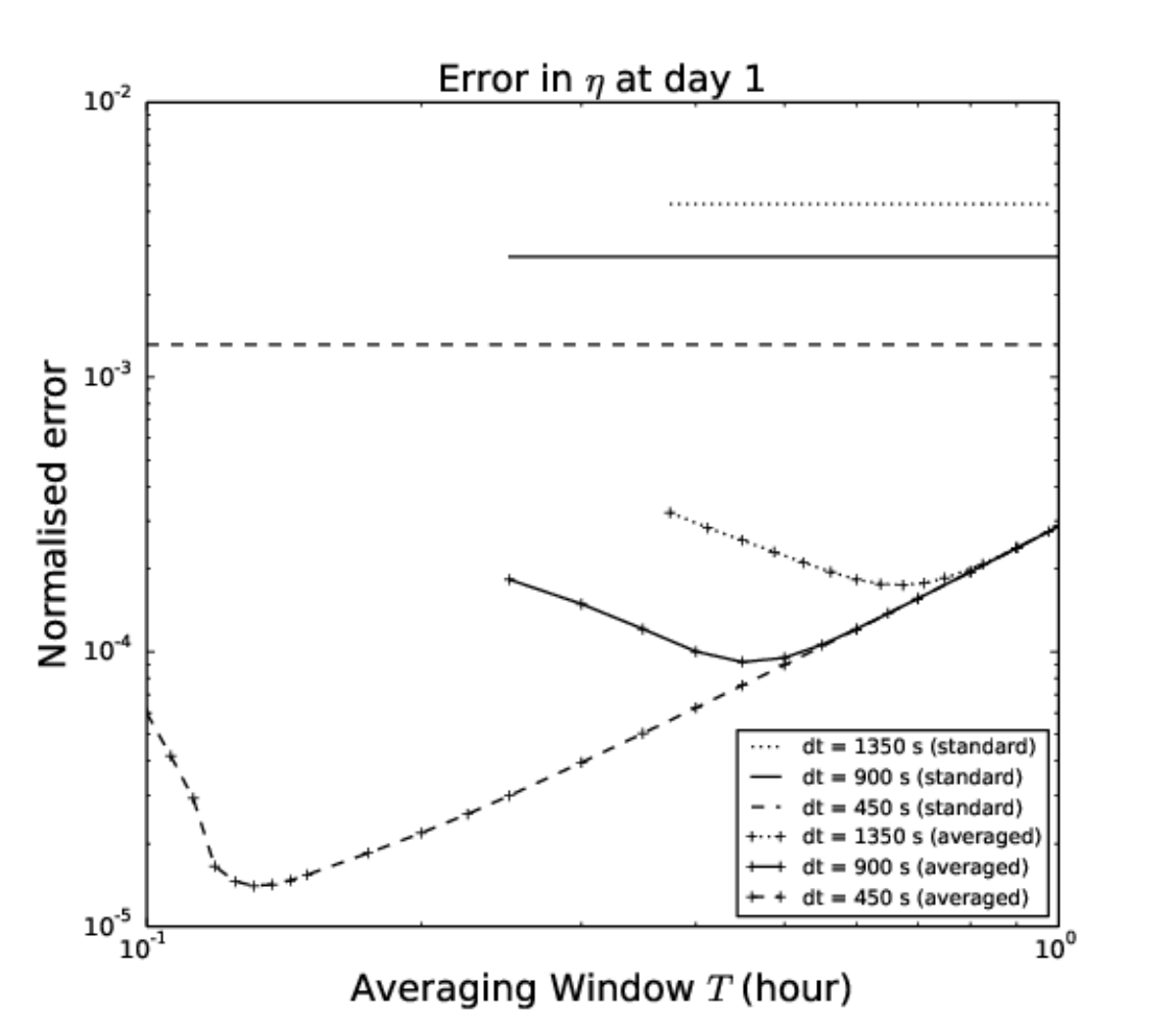} 
  \caption{Normalised errors in the averaged model at day 1 using
    different time step size: (left) the H(div) norm of $\MM{u}$ and
    (right) the $L_2$ norm of $\eta$, normalised by the norms of the
    reference solution. Dashed, solid and dotted {curves} show the
    {errors from the averaged model} using $\Delta t =$ 450, 900 and
    1350 s, respectively.  {Dashed, solid and dotted straight lines show the
      errors from the standard semi-implicit model
      \citep{gibson2019compatible} using $\Delta t =$ 450, 900 and
      1350 s, respectively.}
  \label{fig:convergence_plots}}
\end{figure}

Finally, we examine the accuracy of the averaged model when using different time step size $\Delta t$. {The curves with markers} in Figure \ref{fig:convergence_plots} show the H(div) norm of $\MM{u}$ and the $L_2$ norm of $\eta$ at day 1 {in the averaged model}, normalised by the norms of the reference solution, when using three different time steps: $\Delta t =$ 450, 900 and 1350 s. The averaging window $T$ was changed between $T_{\rm min} \leq T \leq 1$ hour, where $T_{\rm min} = $ 0.1, 0.25 and 0.375 hours for the results using $\Delta t =$ 450, 900 and 1350 s, respectively. When using $T$ smaller than $T_{\rm min}$ that is corresponding to each time step size, the model blows up within 1 day due to the timestepping errors. The results show that, for both variables, the minimum error as well as the optimal averaging window size decrease as the time step is reduced. As the averaging window increases, the amplitudes of the error become almost identical to each other regardless of the time step size. This result confirms that the averaged model is more accurate with a smaller time step when the chosen averaging window is smaller or similar to $\Delta t$, and that the size of $\Delta t$ doesn't affect the accuracy with larger averaging windows where the averaging error dominates.

{
Also shown in Figure \ref{fig:convergence_plots} as lines without markers are the same errors at day 1 in the standard semi-implicit model \citep{gibson2019compatible}, also normalised by the norms of the reference solution, when using the same three different time steps: $\Delta t =$ 450, 900 and 1350 s. As the standard model does not have averaging windows, the errors are shown as straight lines regardless of the size of the averaging window. It is clear that, for all the three time steps used in this test, the solutions from the averaged model is more accurate than those from the standard model using the same time-step size, for the range of the averaging windows shown in Figure \ref{fig:convergence_plots}. In other words, the averaged model would allow us to use a larger time step to achieve the same level of accuracy of the standard model when the averaging windows are chosen near the optimum.
}

\section{Summary and outlook}
\label{sec:conclusion} In this paper we presented a phase averaging framework
for the rotating shallow water equations, and a time integration
methodology for it. The new framework includes overcoming two key
technical challenges for finite element methods on the sphere: the
development of a stable numerical matrix exponential used for the
mappings and a parallel phase averaging procedure. We presented
proof-of-concept results from the rotating shallow water equations and
analysed their errors, which confirm that there is an optimal
averaging window value for a given time step size $\Delta
t$. Naturally, the optimal averaging window for both height and
velocity fields combined depends on the choice of norm. Critically, we
observe that the combined time discretisation and averaging error for
the averaged model is much smaller than the time discretisation error
in a semi-implicit method applied to the same semidiscretisation,
illustrating the benefits of the approach. \add{This is a very
  significant result, because it suggests that phase averaging could
  in itself be used as a time parallel algorithm for the rotating
  shallow water equations on the sphere (and perhaps three dimensional
  models), without necessarily needing corrections through the
  ParaReal algorithm, as proposed in \citet{haut2014asymptotic}.}
  
In future work we will explore the combination of phase averaging
methods with implicit or split timestep methods that allow us to take
even larger timesteps, will incorporate parallel rational
approximation techniques to speed up the exponential evaluations
\citep{haut2016high}, and will undertake parallel performance
benchmarks.

\section*{\large Acknowledgement}
We are grateful for funding from EPSRC under grant EP/R029628/1.
This work used the Isambard 2 UK National Tier-2 HPC Service (\url{http://gw4.ac.uk/isambard/}) operated by GW4 and the UK Met Office, and funded by EPSRC (EP/T022078/1).

\bibliography{averaged_sw}

\begin{thebibliography}{24}
\providecommand{\natexlab}[1]{#1}
\providecommand{\url}[1]{\texttt{#1}}
\providecommand{\urlprefix}{}

\bibitem[{Sanders et~al.(2007)Jan A. Sanders and Ferdinand Verhulst and James
  Murdock}]{Sanders_etal_2007}
Sanders JA, Verhulst F, Murdock J.
\newblock Averaging Methods in Nonlinear Dynamical Systems.
\newblock 2 ed. Springer New York, NY; 2007.

\bibitem[{Majda and Embid(1998)Majda, Andrew J and Embid,
  Pedro}]{majda1998averaging}
Majda AJ, Embid P.
\newblock Averaging over fast gravity waves for geophysical flows with
  unbalanced initial data.
\newblock Theoretical and computational fluid dynamics 1998;11(3):155--169.

\bibitem[{Klainerman and Majda(1981)S. Klainerman and A. J. Majda}]{KlMa1981}
Klainerman S, Majda AJ.
\newblock Singular limits of quasilinear hyperbolic systems with large
  parameters and the incompressible limit of compressible fluids.
\newblock Communications in Pure and Applied Mathematics 1981;34(4):481--524.

\bibitem[{Wagner and Young(2015)Wagner, G. L. and Young, W.
  R.}]{wagner_young_2015}
Wagner GL, Young WR.
\newblock Available potential vorticity and wave-averaged quasi-geostrophic
  flow.
\newblock Journal of Fluid Mechanics 2015;785:401–424.

\bibitem[{Kafiabad et~al.(2021)Kafiabad, Hossein A. and Vanneste, Jacques and
  Young, William R.}]{kafiabad_vanneste_young_2021}
Kafiabad HA, Vanneste J, Young WR.
\newblock Wave-averaged balance: a simple example.
\newblock Journal of Fluid Mechanics 2021;911:R1.

\bibitem[{Haut and Wingate(2014)Haut, Terry and Wingate,
  Beth}]{haut2014asymptotic}
Haut T, Wingate B.
\newblock An asymptotic parallel-in-time method for highly oscillatory {PDE}s.
\newblock SIAM Journal on Scientific Computing 2014;36(2):A693--A713.

\bibitem[{Peddle et~al.(2019)A. Peddle and T. Haut and B.
  Wingate}]{Peddle2019PararealCF}
Peddle A, Haut T, Wingate B.
\newblock Parareal Convergence for Oscillatory {PDE}s with Finite Time-Scale
  Separation.
\newblock SIAM J Sci Comput 2019;41:A3476--A3497.

\bibitem[{Ariel et~al.(2016)Ariel, G. and Kim, S. J. and Tsai,
  R.}]{Ariel_etal_16}
Ariel G, Kim SJ, Tsai R.
\newblock Parareal methods for highly oscillatory dynamical systems.
\newblock SIAM Journal on Scientific Computing 2016;38(6):A3540 --– A3564.

\bibitem[{Bauer et~al.(2022)Bauer, Werner and Cotter, Colin J and Wingate,
  Beth}]{bauer2021higher}
Bauer W, Cotter CJ, Wingate B.
\newblock Higher order phase averaging for highly oscillatory systems.
\newblock To appear in SIAM Journal of Multiscale Methods 2022;.

\bibitem[{Ong et~al.(2016)Ong, Benjamin W and Haynes, Ronald D and Ladd,
  Kyle}]{ong2016algorithm}
Ong BW, Haynes RD, Ladd K.
\newblock Algorithm 965: {RIDC} methods: A family of parallel time integrators.
\newblock ACM Transactions on Mathematical Software (TOMS) 2016;43(1):1--13.

\bibitem[{Minion(2011)Minion, Michael}]{minion2011hybrid}
Minion M.
\newblock A hybrid parareal spectral deferred corrections method.
\newblock Communications in Applied Mathematics and Computational Science
  2011;5(2):265--301.

\bibitem[{Hochbruck and Ostermann(2010)Hochbruck, Marlis and Ostermann,
  Alexander}]{hochbruck2010exponential}
Hochbruck M, Ostermann A.
\newblock Exponential integrators.
\newblock Acta Numerica 2010;19:209--286.

\bibitem[{Pieper et~al.(2019)Pieper, Konstantin and Sockwell, K Chad and
  Gunzburger, Max}]{pieper2019exponential}
Pieper K, Sockwell KC, Gunzburger M.
\newblock Exponential time differencing for mimetic multilayer ocean models.
\newblock Journal of Computational Physics 2019;398:108900.

\bibitem[{Trefethen(2019)Trefethen, Lloyd N}]{trefethen2019approximation}
Trefethen LN.
\newblock Approximation Theory and Approximation Practice, Extended Edition.
\newblock SIAM; 2019.

\bibitem[{Gander and G\"uttel(2013)Gander, Martin J and G\"uttel,
  Stefan}]{gander2013paraexp}
Gander MJ, G\"uttel S.
\newblock {PARAEXP}: A parallel integrator for linear initial-value problems.
\newblock SIAM Journal on Scientific Computing 2013;35(2):C123--C142.

\bibitem[{Haut et~al.(2016)Haut, Terry S and Babb, T and Martinsson, PG and
  Wingate, BA}]{haut2016high}
Haut TS, Babb T, Martinsson P, Wingate B.
\newblock A high-order time-parallel scheme for solving wave propagation
  problems via the direct construction of an approximate time-evolution
  operator.
\newblock IMA Journal of Numerical Analysis 2016;36(2):688--716.

\bibitem[{Schreiber et~al.(2018)Schreiber, Martin and Peixoto, Pedro S and
  Haut, Terry and Wingate, Beth}]{schreiber2018beyond}
Schreiber M, Peixoto PS, Haut T, Wingate B.
\newblock Beyond spatial scalability limitations with a massively parallel
  method for linear oscillatory problems.
\newblock The International Journal of High Performance Computing Applications
  2018;32(6):913--933.

\bibitem[{Croci and Mu{\~n}oz-Matute(2022)Croci, Matteo and Mu{\~n}oz-Matute,
  Judit}]{croci2022exploiting}
Croci M, Mu{\~n}oz-Matute J.
\newblock Exploiting Kronecker structure in exponential integrators: fast
  approximation of the action of $\backslash$$varphi $-functions of matrices
  via quadrature.
\newblock arXiv preprint arXiv:221100696 2022;.

\bibitem[{Gibson et~al.(2019)Gibson, Thomas H and McRae, Andrew TT and Cotter,
  Colin J and Mitchell, Lawrence and Ham, David A}]{gibson2019compatible}
Gibson TH, McRae AT, Cotter CJ, Mitchell L, Ham DA.
\newblock Compatible Finite Element Methods for Geophysical Flows: Automation
  and Implementation Using Firedrake.
\newblock Springer Nature; 2019.

\bibitem[{Rathgeber et~al.(2016)Rathgeber, Florian and Ham, David A and
  Mitchell, Lawrence and Lange, Michael and Luporini, Fabio and McRae, Andrew
  TT and Bercea, Gheorghe-Teodor and Markall, Graham R and Kelly, Paul
  HJ}]{rathgeber2016firedrake}
Rathgeber F, Ham DA, Mitchell L, Lange M, Luporini F, McRae AT, et~al.
\newblock Firedrake: automating the finite element method by composing
  abstractions.
\newblock ACM Transactions on Mathematical Software (TOMS) 2016;43(3):1--27.

\bibitem[{Balay et~al.(2021)Balay, Satish and Abhyankar, Shrirang and Adams,
  Mark and Brown, Jed and Brune, Peter and Buschelman, Kris and Dalcin,
  Lisandro and Dener, Alp and Eijkhout, Victor and Gropp, W and
  others}]{balay2021petsc}
Balay S, Abhyankar S, Adams M, Brown J, Brune P, Buschelman K, et~al.
\newblock PETSc users manual: revision 3.15 2021;.

\bibitem[{Williamson et~al.(1992)Williamson, David L and Drake, John B and
  Hack, James J and Jakob, R{\"u}diger and Swarztrauber, Paul
  N}]{williamson1992standard}
Williamson DL, Drake JB, Hack JJ, Jakob R, Swarztrauber PN.
\newblock A standard test set for numerical approximations to the shallow water
  equations in spherical geometry.
\newblock Journal of Computational Physics 1992;102(1):211--224.

\bibitem[{Thuburn et~al.(2014)Thuburn, J and Cotter, CJ and Dubos,
  T}]{thuburn2014mimetic}
Thuburn J, Cotter C, Dubos T.
\newblock A mimetic, semi-implicit, forward-in-time, finite volume shallow
  water model: comparison of hexagonal--icosahedral and cubed-sphere grids.
\newblock Geoscientific Model Development 2014;7(3):909--929.

\bibitem[{Shipton et~al.(2018)Shipton, Jemma and Gibson, Thomas H and Cotter,
  Colin J}]{shipton2018higher}
Shipton J, Gibson TH, Cotter CJ.
\newblock Higher-order compatible finite element schemes for the nonlinear
  rotating shallow water equations on the sphere.
\newblock Journal of Computational Physics 2018;375:1121--1137.

\end{thebibliography}

\end{document}